Revised December 2016

# A Note on Graded Rings and Modules

*by* Nicholas Phat Nguyen

E-mail address: nicholas.pn@gmail.com

**Abstract.** In this note, we consider a situation that is generally used as an intermediate technical step in proving the Artin-Rees lemma but otherwise is not much discussed in introductory accounts of commutative algebra. I hope to show in this note that such technical step deserves more recognition and emphasis in any introduction to commutative algebra because it can be used to prove some significant results in a straight-forward manner, including a generalization of a theorem by Pierre Cartier and John Tate.

A.  INTRODUCTION:

We begin by defining our basic objects. Let $S = \bigoplus_{k \geq 0} S_k$ be a commutative graded ring, with the usual grading rule $S_i.S_j \subset S_{i+j}$. With this grading rule, $S_0$ is a subring of $S$, and each component $S_k$ can be regarded as a module over $S_0$ under the multiplication in $S$.

In addition, let $M = \bigoplus_{k \geq 0} M_k$ be a direct sum of commutative groups. Assume that we also have a pairing $S \times M \to M$ that makes $M$ into an $S$-module. We say that $M$ is a graded $S$-module if $S_i.M_j \subset M_{i+j}$ for any indices $i$ and $j$. (Obviously, $S$ itself is a graded module over $S$.) An element $x$ of $M_k$ is called homogeneous of degree $k$. With the grading rule, each $M_k$ can naturally be regarded as a module over the subring $S_0$. We say the $S$-grading in $M$ is simple if for all sufficiently large indices $k$, we have $S_1 M_k = M_{k+1}$.

A submodule N of M is a graded module over S if and only if $N = \bigoplus_{k \geq 0} N_k$ where $N_k \subset M_k$ for each index k, i.e., if each homogenous component of an element of N also belongs to N. If a submodule N does not satisfy this grading condition, we will refer to it as an ungraded submodule.

We will assume throughout our discussion the following two conditions:

(1) As a ring, S is generated by $S_0$ and $S_1$;
(2) Each $M_k$ is finitely generated as a module over the subring $S_0$.

Condition (1) means that each element of $S_k$ can be written as a sum of products of k homogeneous elements of degree 1, i.e. $S_k = (S_1)^k$. If $S_0$ is a Noetherian ring and $S_1$ is a finitely generated module over $S_0$, then condition (1) implies that S is also a Noetherian ring by the Hilbert's basis theorem. However, unless specifically indicated otherwise, we do not assume in general that $S_0$ is Noetherian.

B.   A TECHNICAL LEMMA:

Most expositions of commutative algebra mention the following technical lemma as a step in proving the Artin-Rees lemma.

*Technical Lemma*: Let S and M be given and satisfy conditions (1) and (2) above, then M is a finitely generated module over S if and only if the S-grading on M is simple.

Other than helping to prove the Artin-Rees lemma, the above technical lemma does not seem to receive any attention in most introductory accounts of commutative algebra. That is rather unfortunate because this technical lemma is conceptually more fundamental and significant than the Artin-Rees lemma per se and can be used to prove a number of other interesting and important results, as I will try to show in this note. [1]

---

[1] Professor Brian Conrad (Stanford University) pointed out to the author that graded modules over graded rings arise in diverse situations in mathematics, such as in modern projective algebraic geometry (e.g., coherent sheaves on projective space) and the combinatorial side of representation theory.   In those

The proof of the above technical lemma is quite simple. For the convenience of the reader, we sketch the proof of that lemma in order to highlight the essential ideas.

Proof of Lemma: Assume first that M is a finitely generated module over S. Then M can be generated by the subgroup $M_0 + M_1 + ... + M_k$ for some index k. That implies that for any index $j \geq k$, we have $M_{j+1} = S_{j+1} M_0 + S_j M_1 + ... + S_1 M_j$. Because of condition (1), that can be written as $M_{j+1} = (S_1)^{j+1} M_0 + (S_1)^j M_1 + ... + S_1 M_j \subset S_1 M_j$. Since clearly $S_1 M_j \subset M_{j+1}$ because of the grading condition, we have $S_1 M_j = M_{j+1}$ for any index $j \geq k$. Therefore the S grading on M is simple.

Conversely, assume that the S grading on M is simple, and we have $S_1 M_j = M_{j+1}$ for any index $j \geq k$. Then it is clear that M can be generated as a module over S by the subgroup $M_0 + M_1 + ... + M_k$. Because of condition (2), we can select a finite set of generators from that subgroup that would generate M over S.[2]

As an immediate application of this technical lemma, we can use it to provide a simple proof of the following general version of Krull's intersection theorem.

*Krull's Intersection Theorem*: Let R be a Noetherian ring, **a** an ideal of R, and M a finitely generated R module. If N is any submodule of $\bigcap_{i \geq 0} \mathbf{a}^i M$, then we have $N = \mathbf{a}N$.

Proof: Let $S = R \oplus \mathbf{a}R \oplus \mathbf{a}^2 R \oplus ...$, and let $P = M \oplus \mathbf{a}M \oplus \mathbf{a}^2 M \oplus ...$ S is naturally a graded ring and P is naturally a graded S-module. Moreover, S and P meet conditions (1) and (2). Because the S grading on P is simple by design, P is finitely generated over S. The Noetherian hypothesis on R implies that S is Noetherian, and therefore P is also Noetherian. Accordingly, any S-submodule of P is also finitely generated over S.

---

situations, it is natural to expect the technical lemma (in one form or another) to be a basic tool. However, in introductory accounts of commutative algebra, that significance of the technical lemma is not easily seen.

[2] The above outline actually shows that the technical lemma has a more precise version: with condition (1), M is finitely generated over S implies that M is simply graded; and with condition (2), M is simply graded implies M is finitely generated over S.

Let $Q = N \oplus N \oplus \ldots$ be the graded S module where each component is equal to N. Because N is contained in $\mathbf{a}^i M$ for any index i, Q is a graded submodule of P, and hence it is finitely generated over S. By the technical lemma, the S grading on N is simple, which means that $\mathbf{a}N = N$.

The above proof is conceptually more straight-forward than the usual proof of Krull's intersection theorem using the Artin-Rees lemma. In fact, if we unpack the technical terms, the above proof is essentially just a more abstract version of the proof by Hervé Perdry in [Perdry 2004].

C.  A FINITENESS THEOREM:

We will apply the above technical lemma to prove the following generalization of a theorem by Cartier and Tate from [Cartier Tate 1978].

As before, let M be a graded module over the graded ring S that satisfies conditions (1) and (2). Given a submodule N of M (whether or not graded), we will sometimes want to know if the quotient module M/N is finitely generated over the subring $S_0$.

For convenience, we will refer to a graded M module as short if the components $M_k$ are zero for all sufficiently large indices k. In contrast, a graded module M is called long if we can find arbitrarily large index j such that $M_j$ is non-zero. Clearly a graded module M is finitely generated over the subring $S_0$ if and only if it is a short module, in light of condition (2).

If a simply graded module has a zero component for a sufficiently large index when the grading has become simple, then all components are obviously zero from that point on. Therefore in the case of simple grading, a graded module is long if and only if all of its components are non-zero beyond a sufficiently large index.

We will refer to any S-submodule N of M (whether or not graded) as saturated if N contains all the components $M_k + M_{k+1} + \ldots$ for some index k. For any saturated submodule N of M (whether or not graded), the quotient module M/N is clearly a finitely generated module over the subring $S_0$ by virtue of condition (2) above. If the submodule N is graded, then the converse is also true: if the quotient graded module M/N is finitely generated over the subring $S_0$, then M/N must be a short graded module, and therefore N is a saturated submodule of M.

Given a long graded module M, we can ask whether there is a non-saturated submodule N such that the quotient module M/N is finitely generated over the subring $S_0$. Such a non-saturated submodule, if it exists, must necessary be non-graded because a graded non-saturated submodule N will give us a graded quotient module M/N that is also long and hence cannot be finitely generated over the subring $S_0$.

*Theorem*. Let M be a long graded module over S. Assume that:
  (a) S and M meet conditions (1) and (2) above;
  (b) $S_1$ is a finitely generated module over $S_0$; and
  (c) the S grading on M is simple.
In that case, we can always find a non-saturated S-submodule L of M such that the quotient S-module M/L is finitely generated as a module over the subring $S_0$.

*Example*: Let $S = F[X_0, X_1, ..., X_n]$ be the polynomial ring in variables $X_0, X_1, ..., X_n$ over a field F, with the natural grading given by homogeneous degree. Given a homogeneous ideal J of S that is non-saturated, the graded quotient ring S/J is a long graded S-module that meets the conditions of the theorem. The theorem tells us that there is a non-saturated ideal K containing J such that the quotient ring S/K is finitely generated as an F-module, i.e., a finite-dimensional F-vector space. Accordingly, we can find a surjective ring homomorphism of S/J onto a finite-dimensional algebra over F such that the images of $X_0$, $X_1, ..., X_n$ are not all nilpotent. Because a finite-dimensional F-algebra is Artinian, any prime ideal is also maximal and therefore we can always choose a maximal ideal that avoids a given non-nilpotent element. Accordingly, we can factor out a suitable maximal ideal of such finite dimensional F-algebra to obtain a ring homomorphism of S/J into an algebraic extension of F where the images of $X_0, X_1, ..., X_n$ are not all zero. Let $x_0, x_1, ..., x_n$ be such images. For any polynomial f in J, we have $f(x_0, x_1, ..., x_n) = 0$. What we have then is the homogeneous version of Hilbert's zeros theorem (Hilbert's Nullstellensatz): any non-saturated homogeneous ideal J of S has a zero in the projective space $\mathbf{P}_n$ defined over an algebraic closure of F.

In a 1978 paper, Pierre Cartier and John Tate gave a conceptually simple proof of the homogeneous version of Hilbert's Nullstellensatz and the elimination theorem in algebraic geometry using the graded ring structure of $F[X_0, X_1, ..., X_n]$ and its quotient rings. See [Cartier and Tate 1978]. What we are doing in this note is generalizing their arguments and results to more general graded modules. A key step in the proof will be the technical lemma.

Proof of Theorem. We outline a proof of the above finiteness theorem in a series of steps.

i. By hypothesis the graded module M is long and also simply graded. That means each of its components beyond a sufficiently large index must be non-zero.

ii. Take a generating set A of the module $S_1$ over $S_0$. By hypothesis, we can choose A to be finite.

iii. A is non-empty because in that case $S_1 = 0$, $S = S_0$ and all components of M are zero beyond a certain index because of the simple grading condition, contrary to the assumptions of the theorem.

iv. For each subset B of A, let N be the subgroup of M defined by letting $N_0 = 0$ and $N_{k+1} = BM_k$ for any $k \geq 0$. N is a graded S-submodule of M. Moreover, the grading on N is also simple because the grading on M is simple.

v. N obviously depends on the chosen subset B. If B is the empty set, then $N = 0$ and $M/N = M$ is still a long graded module. At the other end, if $B = A$, then for all sufficiently large indices j, we claim $N_j = M_j$. Although the initial components of N and M may be different, we must have $N_{k+1} = AM_k = S_1M_k = M_{k+1}$ at any index k where M has become simply graded. Accordingly, M/N is a short graded module when $B = A$.

vi. For each subset B of A, we have a dichotomy for the corresponding submodule N: either $N_k = M_k$ for all sufficiently large indices k (as for example when $B = A$), or $N_k$ is strictly smaller than $M_k$ for all sufficiently large indices k (as for example when B is the empty set). Indeed, if $N_{k+1} = BM_k = M_{k+1}$ at an index k where M

      has become simply graded, then $N_{k+2} = BM_{k+1} = B(S_1 M_k) = S_1(BM_k) = S_1 M_{k+1} = M_{k+2}$. By induction, both $N_j$ and $M_j$ coincide for any index $j > k$.

vii. Because A is finite, we can chose a proper subset B between $\emptyset$ and A with the property that the graded quotient S-module M/N is long (i.e. has no zero component beyond a sufficiently large index), but if we add any element x in A but $\notin$ B to B, then the corresponding graded quotient module M/N is short (i.e., has all its components vanishing for sufficiently large indices.)

viii. It is enough for us to prove our theorem for the graded S-module M/N constructed with such a maximal subset B. That graded S-module is still long and continues to satisfy all the hypothesis in the theorem. Any non-saturated submodule of M/N that has the desired properties will give us by pull-back a non-saturated S-submodule of M with the same desired properties.

ix. The maximal character of N means that for any element x in A but $\notin$ B, $x(M/N)_k = (M/N)_{k+1}$ for all sufficiently large indices k, i.e., $M_{k+1} = xM_k + N_{k+1} = xM_k + BM_k$. But that follows from our choice of B because when we use $B \cup (x)$ as the generating set, the resulting graded submodule has the same components as M for all sufficiently large indices k.

x. M/N is a graded module over S. It can also be regarded as a graded module over $S_0[x]$, where $S_0[x]$ is the graded subring of S generated by $S_0$ and x. Because $x(M/N)_k = (M/N)_{k+1}$ for all sufficiently large indices k, the grading of $S_0[x]$ on M/N is simple.

xi. By the technical lemma, that means M/N is finitely generated as a module over the subring $S_0[x]$.

xii. To simplify notation, write $P = M/N$. We claim that the (non-graded) submodule $(1 - x)P$ of P has the desired properties. That means (i) the S-submodule $(1-x)P$ of P is non-saturated, and (ii) the quotient S-module $P/(1 - x)P$ is finitely generated as a module over $S_0$.

xiii. To say the S-submodule $(1-x)P$ is non-saturated in P means that for any index k no matter how large, we can pick an index $j \geq k$ such that the component $P_j$ of P is not contained entirely in $(1-x)P$. We can assume that any such $P_j$ is non-zero

(P is a long graded module) and moreover $xP_j = P_{j+1}$ if k is large enough (from the choice of maximal B).

xiv. If a homogeneous element **v** of $P_j$ is equal to $(1-x)\mathbf{u}$ for some **u** of P, then **u** must be equal to $\mathbf{v} + x\mathbf{v} + x^2\mathbf{v} + \ldots + x^n\mathbf{v}$ with $x^{n+1}\mathbf{v} = 0$ for some natural integer n, as we can readily check by expressing u in terms of its homogenous components. Accordingly, if $P_j$ is contained entirely in $(1-x)P$, then each non-zero element **v** of $P_j$ can be annihilated by a power of x. Because $P_j$ is finitely generated over $S_0$ (condition (2) applied to the module $P = M/N$), that means for large enough power n of x, $x^n P_j = 0$, and therefore all components $P_j$ are zero beyond a sufficiently large index, contrary to our assumption. This argument shows that $(1-x)P$ is a non-saturated S-submodule of P.

xv. To simplify notation, write R for the graded subring $S_0[x]$ of S. Note that the quotient S-module $P/(1-x)P$ is naturally isomorphic (as S module) to the tensor product $P \otimes_R R/(1-x)R$ where the tensoring is relative to the R-module structure of P and $R/(1-x)R$. Indeed, if we start with the following short exact sequence of R-modules:

$$0 \to (1-x)R \to R \to R/(1-x)R \to 0,$$

and then tensor (over R) with the module P (viewed as both an R module and an S module), we obtain the following right exact sequence of S-modules:

$$P \otimes_R (1-x)R \to P \to P \otimes_R R/(1-x)R \to 0.$$

Because the image of $P \otimes_R (1-x)R \to P$ is just the submodule $(1-x)P$ of P, we have a natural isomorphism between $P/(1-x)P$ and $P \otimes_R R/(1-x)R$.

xvi. We can now readily show that the module $P \otimes_R R/(1-x)R$ is finitely generated as a module over the subring $S_0$. Note that P is finitely generated over R as we discovered above. So it is sufficient to show that the R module $R/(1-x)R$ is finitely generated as a module over the subring $S_0$ of $R = S_0[x]$.

xvii. For any polynomial f in $R = S_0[x]$, $f \equiv f(1) \mod (1-x)$, and so the natural mapping from $S_0$ to $R/(1-x)R$ given by $s \to s \mod (1-x)$ is surjective. That means if we regard $R/(1-x)R$ as an $S_0$-module, it is naturally isomorphic to a quotient of $S_0$ itself, and generated by the element 1 mod (1 – x) as a module over $S_0$.

With this last step, the proof of the above finiteness theorem is complete.

**Acknowledgment.** The author would like to thank Professor Brian Conrad (Stanford University) for his helpful and valuable comments in reviewing an earlier draft of the note. The author would also like to express his sincere thanks to Vincent Beck and to Jan R. Strooker for calling his attention to certain errors and typos in earlier versions of this paper.

______________________________________________.

REFERENCES:

[Cartier Tate 1978] P. Cartier and J. Tate, *A simple proof of the main theorem of elimination theory in algebraic geometry*, L'enseignement mathematique, Vol. 24 at p. 311 (1978).

[Milne 2013] J.S. Milne, A *Primer of Commutative Algebra* (May 2013 version), available at www.jmilne.org/math.

[Perdry 2004] Hervé Perdry. *An elementary proof of Krull's intersection theorem*, The American Mathematical Monthly, 111(4):356-357, 2004.